\documentclass[10pt,twocolumn]{amsart}
\usepackage[T1]{fontenc}
\usepackage[utf8]{inputenc}
\usepackage{amsmath,amssymb,amsthm,mathtools,bm}

\usepackage{newtxtext,newtxmath}
\usepackage[letterpaper,margin=0.66in]{geometry}
\usepackage{microtype,hyperref}
\usepackage{multicol}
\usepackage{tikz}
\hypersetup{colorlinks=true,linkcolor=blue,citecolor=blue,urlcolor=blue}
\newcommand{\R}{\mathbb R}
\newcommand{\E}{\mathbb E}
\newcommand{\Ph}{\mathbf P_h}
\newcommand{\Ps}{\mathbf P_\sigma}
\newcommand{\PU}{\mathbf P_U}
\newcommand{\EU}{\mathbf E_U}
\newcommand{\Qn}{\mathbf Q_n^{(\rho)}}
\newcommand{\PPP}{\mathbf{PPP}}
\newcommand{\1}[1]{\bm 1_{\{#1\}}}
\DeclareMathOperator{\sign}{sign}
\DeclareMathOperator{\var}{Var}

\title{REM universality and Poisson--Dirichlet Gibbs weights\vspace{-1mm}\\ for linear random energy}

\author{Francesco Concetti\textsuperscript{1,*}, Simone Franchini\textsuperscript{2}}

\address{
\textnormal{\textsuperscript{1}~Faculty of Mathematics and Computer Science, UniDistance Suisse, 3900 Brig, Switzerland}\\
\textnormal{\textsuperscript{2}~CNR-ISTC, Via Gian Domenico Romagnosi 18, 00196 Rome, Italy}\\
\textnormal{\textsuperscript{*}~Contact: \href{mailto:francesco.concetti@unidistance.ch}{francesco.concetti@unidistance.ch}}
}

\date{}

\begin{document}
\twocolumn[
\begin{@twocolumnfalse}
\begin{abstract}
We study the Hamiltonian $H_n(h,\sigma)=\sum_{i=1}^n h_i(\sigma_i-m),
$ where $(h_i)$ are i.i.d.\ real random variables and $(\sigma_i)$ are i.i.d.\ Ising spins. We consider the energy levels obtained after an independent thinning that retains an exponential number of configurations ($e^{O(n)}$). We prove that, after an $(h_i)$-dependent centering, the resulting point process converges in distribution to a Poisson point process with exponential intensity. Thus, the energy levels asymptotically has the one of the Random Energy Model (REM). Our results extend previous ones, where REM universality for this model was established only either for energy fluctuations of order $e^{-O(n)}$ or for $e^{o(\sqrt n)}$ randomly selected configurations. We also identify the limiting Gibbs weights, which converge to a Poisson--Dirichlet law, and the quenched free energy, which exhibits a freezing transition at $\beta=\tilde\lambda$. The proofs are presented here in compressed form; full details are given in the companion preprint.
\end{abstract}
\maketitle

\vspace{2mm}
\end{@twocolumnfalse}
]
\section*{Introduction}

REM universality asks whether, after a suitable centering, the energy levels of a disordered Hamiltonian on $n$ spins exhibit, as $n\to\infty$, the same local statistics as Derrida's Random Energy Model (REM) \cite{Derrida}. For number partitioning and related mean-field Hamiltonians, Borgs, Chayes, Mertens, and Nair, and independently Bovier and Kurkova, showed that the REM picture emerges inside an energy window shrinking exponentially fast to zero with the system size \cite{Borgs,BovierK}. Ben Arous, Gayrard, and Kuptsov later introduced \emph{REM universality by dilution}, proving asymptotic REM behavior for order-one fluctuations of $e^{o(\sqrt n)}$ randomly selected energy levels \cite{Arous}. REM universality without dilution and without shrinking windows was subsequently obtained by Arguin and Kistler \cite{ArguinKistler} for a model given by the sum of a linear random Hamiltonian and an independent REM Hamiltonian.

In the present work, we establish REM universality for order-one fluctuations among $e^{O(n)}$ randomly selected energy levels of a purely linear random Hamiltonian, with no added REM term.

The importance of REM universality stems from the fact that, in the REM, the properly centered extremal energy levels converge to a Poisson point process with exponential intensity measure, yielding a canonical and tractable limit for disordered spectra. In the low-temperature phase, the distribution of the ranked Gibbs weights of the REM converges to the one-parameter Poisson--Dirichlet law $\mathrm{PD}(\theta,0)$; see \cite[Definition~(1), Proposition~(10)]{PitmanYor}.

In spin-glass language, this is precisely the law of the ranked pure-state weights in a one-step replica-symmetry-breaking ($1$RSB) Ruelle probability cascade, that is, the simplest rigorous realization of the $1$RSB Parisi picture \cite{Ruelle,ArguinRPC,As2}.  The connection between REM universality and the Parisi solution has recently attracted renewed attention in the physics literature on simplified Parisi descriptions and mean-field spin glasses, motivated by a conjectural REM-universality scenario \cite{Fra21,Fra23,Fra25}.

The REM has also recently emerged as a powerful tool in the statistical-mechanics analysis of retrieval thresholds and capacity bounds for a class of Hopfield-type models with exponential storage capacity, namely dense associative memories \cite{Hopfield,LucibelloMezard}.

This letter gives a condensed presentation of the rigorous results established in our companion preprint \cite{CF26}. We study the purely linear random Hamiltonian
\begin{equation}
H_n(h,\sigma):=\sum_{i=1}^n h_i(\sigma_i-m),
\label{eq:Hn}
\end{equation}
where $m\in(-1,1)$ is fixed, $h=(h_i)_{i\ge1}$ are i.i.d.\ real random variables, and $\sigma=(\sigma_i)_{i\ge1}$ are independent $\{-1,1\}$-valued spins, independent of $h$, with
\[
\Ps(\sigma_1=1)=\frac{1+m}{2},
\qquad
\Ps(\sigma_1=-1)=\frac{1-m}{2}.
\]
For each $n$, the family $\{H_n(h,\sigma):\sigma\in\{-1,1\}^n\}$ is strongly correlated, since all energies are linear functionals of the same random environment. Following the thinning perspective of \cite{Arous}, we retain only a random subset of configurations, but in contrast with \cite{Arous} the thinning keeps an exponential number of them ($e^{O(n)}$ configurations).

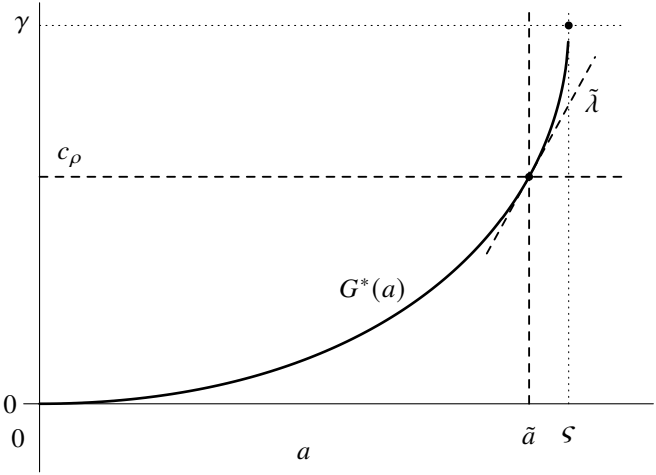
\begin{figure}[htbp]
\centering
\begin{tikzpicture}[x=7cm,y=5cm,line cap=round,line join=round]

% Parameters
\pgfmathsetmacro{\varsig}{1.00}
\pgfmathsetmacro{\gammaVal}{1.00}
\pgfmathsetmacro{\cRho}{0.60}
\pgfmathsetmacro{\atilde}{0.92519}
\pgfmathsetmacro{\lambdaTilde}{2.53727}
\pgfmathsetmacro{\tangentLeft}{\atilde - 0.08}
\pgfmathsetmacro{\tangentRight}{1.05}
\pgfmathsetmacro{\slopeLabelX}{1.03}
\pgfmathsetmacro{\slopeLabelY}{\cRho + \lambdaTilde*(\slopeLabelX-\atilde) - 0.10}

% Axes
\draw[black] (-0.03,0) -- (1.16,0);
\draw[black] (0,-0.18) -- (0,1.06);

% Main curve: G*(a) = 1 - 3/2 sqrt(1-a) + 1/2 (1-a)^{3/2}
% It satisfies G*(0)=0, (G*)'(0)=0, G*(1)=1, and (G*)'(1^-)=+infty.
\draw[black, line width=1pt, domain=0:1, samples=300, smooth]
  plot (\x, {1 - 1.5*sqrt(1-\x) + 0.5*(1-\x)^(1.5)});

% Horizontal/vertical reference lines
\draw[black, dashed,line width=0.7pt] (0,\cRho) -- (1.10,\cRho);
\draw[black, dashed,line width=0.7pt] (\atilde,0) -- (\atilde,1.04);
\draw[black, dotted] (\varsig,0) -- (\varsig,1.04);
\draw[black, dotted] (0,\gammaVal) -- (1.10,\gammaVal);

% Key points
%\fill[black] (0,0) circle (1.5pt);
\fill[black] (\atilde,\cRho) circle (1.5pt);
\fill[black] (\varsig,\gammaVal) circle (1.5pt);

% Tangent line at the intersection point
\draw[black, dashed,line width=0.7pt]
  (\tangentLeft,{\cRho + \lambdaTilde*(\tangentLeft-\atilde)}) --
  (\tangentRight,{\cRho + \lambdaTilde*(\tangentRight-\atilde)});

% Labels (kept away from the lines)
\node at (0.63,0.30) {$G^*(a)$};
\node[above left] at (0.10,\cRho) {$c_{\rho}$};
\node[left] at (0,\gammaVal) {$\gamma$};
\node[below=13pt] at (0.50,0) {$a$};
\node[below=5pt] at (\atilde,0) {$\tilde a$};
\node[below=5pt] at (\varsig,0) {$\varsigma$};
\node[left=5pt] at (0,0) {$0$};
\node[below=13pt,left=2pt] at (0,0) {$0$};
\node[below=-5pt,right=-5pt] at (\slopeLabelX,\slopeLabelY) {$\,\tilde\lambda$};

\end{tikzpicture}
\caption{The curve $G^*$ denotes the quenched large-deviation rate function of $H_n$ conditional on the environment $h$. The parameter $\tilde a$ is defined by $G^*(\tilde a)=c_\rho$, and $\tilde\lambda$ is the slope of $G^*$ at $\tilde a$.}\label{fig:1}
\end{figure}

More precisely, for $\rho\in(0,1)$ we set
\[
c_\rho:=\rho\bigl(\log 2-\log(1+|m|)\bigr),
\]
and for each $\tau\in\{-1,1\}^n$ we define
\begin{equation}
\Qn(\tau):=e^{nc_\rho}\Ps(\tau)
=
e^{nc_\rho}\prod_{i=1}^n \frac{1+\tau_i m}{2}.
\label{eq:Qn}
\end{equation}
Let $(U_\sigma)_{\sigma\in\{-1,1\}^n}$ be i.i.d. uniform on $[0,1]$, independent of everything else. The retained configurations are
\[
\mathcal G_n(U):=\bigl\{\sigma\in\{-1,1\}^n:U_\sigma<\Qn(\sigma)\bigr\}.
\]
so that the average number of retained configurations is
\[
\E[|\mathcal G_n(U)|]=e^{nc_\rho}.
\]
Assuming the regularity hypotheses of Assumption~1.1 in \cite{CF26} on the law of $h_1$---in particular, the linear small-ball condition, a nontrivial absolutely continuous component with density bounded from below on some interval, and $\E|h_1|^3<\infty$---we prove that there exists an $h$-measurable centering $A_n(h)$ such that the thinned point process
\begin{equation}
\mathbf H_n:=\sum_{\sigma\in\mathcal G_n(U)}
\delta_{H_n(h,\sigma)-A_n(h)}
\label{eq:Hpoint}
\end{equation}
converges in distribution to a Poisson point process on $\R$ with exponential intensity measure
\begin{equation}
\mathcal D_{\tilde\lambda}(dx):=e^{-\tilde\lambda x}\,dx.
\label{eq:exponential}
\end{equation}
A graphical characterization of the parameter $\tilde\lambda$ is given in Fig.~\ref{fig:1}.

For $\beta>0$, we define the partition function of the thinned system
\[
Z_n(\beta):=\sum_{\sigma\in\mathcal G_n(U)}e^{\beta H_n(h,\sigma)},
\]
that yields the Gibbs weights
\[
G_n(\sigma):=\frac{e^{\beta H_n(h,\sigma)}}{Z_n(\beta)},
\qquad \sigma\in\mathcal G_n(U).
\]
and the free energy 
\[
F(\beta):=\lim_{n\to \infty}\frac{1}{n}\log Z_n(\beta).
\]
We show that, in the low-temperature regime $\beta>\tilde\lambda$, the convergence of $\mathbf H_n$ to a Poisson point process implies the convergence of the joint distribution of the ranked Gibbs weights to the Poisson--Dirichlet law $\mathrm{PD}(\tilde\lambda/\beta,0)$. In spin-glass language, the Gibbs measure therefore converges to a one-level $1$RSB Ruelle probability cascade. We also compute the limiting free energy explicitly, obtaining a phase transition at $\beta=\tilde\lambda$.

The next sections keep only the backbone of the argument, while the fully rigorous finite-$n$ analysis is given in the companion preprint \cite{CF26}.
\section*{The conditional kernel}
In this section, we analyze the quenched law of a single energy level, namely the conditional distribution of \(H_n(h,\sigma)\) given the environment \(h\).

We introduce the quenched large-deviation rate function that determines both
the centering and the limiting intensity. Set
\[
G(\lambda):=\E\left[\log\!\left(\frac{1+m}{2}e^{\lambda h_1(1-m)}+\frac{1-m}{2}e^{-\lambda h_1(1+m)}\right) \right]\\
\]
and
\begin{equation}
G^*(a):=\sup_{\lambda\in\R}\{\lambda a-G(\lambda)\},
%\label{eq:Gdef}
\end{equation}
We also define the parameters
\[
\varsigma:=-m\E[h_1]+\E[|h_1|],
\]
\[
\gamma:=-\E\log\bigl(1+m\,\sign(h_1)\bigr)+\log 2.
\]
For each \(c\in(0,\gamma)\), there is a unique pair
\((\tilde a,\tilde\lambda)\in(0,\varsigma)\times(0,\infty)\) characterized by
\begin{equation}
G^*(\tilde a)=c,
\qquad
G'(\tilde\lambda)=\tilde a.
\label{eq:solac}
\end{equation}
 Set
\[
\Gamma_n(h):=-\sum_{i=1}^n\log\bigl(1+m\,\sign(h_i)\bigr)+n\log 2.
\]
and consider any \(h\)-measurable sequence \(C_n(h)\) such
that
\[
0<C_n(h)<\Gamma_n(h)\quad\text{eventually},
\qquad
\frac{C_n(h)}{n}\to c
\quad \Ph\text{-a.s.}
\]
The key point is that one can choose an \(h\)-measurable centering \(A_n(h)\) so
that the rescaled laws
\begin{equation}
K_n(h,dx):=e^{C_n(h)}\Ps\bigl(H_n(h,\sigma)-A_n(h)\in [x,x+dx]\bigr)
\label{eq:Kn}
\end{equation}
converge vaguely, \(\Ph\)-almost surely, to the exponential measure
\(\mathcal D_{\tilde\lambda}\). In the next section we will apply this with
\(C_n(h)=n c_\rho\).

To see why this, set
\[
M_n(h,\lambda):=\log\,\E_\sigma\!\left[e^{\lambda H_n(h,\sigma)}\right]
=\sum_{i=1}^n g(\lambda h_i),
\]
and let \(M_n^*(h,\cdot)\) be its Legendre transform. Denoting by
\(\Lambda_n(h,a)\) the unique positive solution of
\(M_n'(h,\Lambda)=a\), we take \(A_n(h)\) so that
\begin{equation}
M_n^*\bigl(h,A_n\bigr)+\frac12\log\Bigl(2\pi M_n''\bigl(h,\Lambda_n(h,A_n)\bigr)\Bigr)
=C_n(h).
\label{eq:Anchoice}
\end{equation}
The existence of such a centering, together with its asymptotic behavior, is
proved rigorously in \cite{CF26}.

The conditional strong large-deviation theorem of Bovier and Mayer
\cite[Theorem~1.6]{BovierMayer} then gives, for every compact \(B\subset\R\),
uniformly in \(x\in B\),
\begin{equation}
\Ps\bigl(H_n\ge A_n+x\bigr)
=
\frac{e^{-M_n^*(h,A_n+x)}}{\Lambda_n^x\sqrt{2\pi M_n''(h,\Lambda_n^x)}}
\bigl(1+o_K(1)\bigr),
\label{eq:sldp}
\end{equation}
where \(\Lambda_n^x=\Lambda_n(h,A_n+x)\).

Since \(x\) stays in an \(O(1)\) window around \(A_n(h)\), Taylor expansion gives
\begin{equation}
M_n^*(h,A_n+x)=M_n^*(h,A_n)+\Lambda_n(h,A_n)x+O(n^{-1}),
\label{eq:taylorMstar0}
\end{equation}
and
\begin{equation}
\Lambda_n(h,A_n+x)=\Lambda_n(h,A_n)+O(n^{-1}),
\label{eq:taylorMstar}
\end{equation}
uniformly on compact sets.

Moreover, \cite{CF26} proves that, \(\Ph\)-almost surely,
\begin{equation}
\label{eq:limits}
\left(\frac{M_n^*(h,A_n)}{n},\frac{A_n(h)}{n},\Lambda_n(h,A_n)\right)
\xrightarrow{\Ph\text{-a.s.}}
\bigl(G^*(\tilde a),\tilde a,\tilde\lambda\bigr).
\end{equation}
Combining \eqref{eq:Anchoice}, \eqref{eq:sldp}, \eqref{eq:taylorMstar0}, and
\eqref{eq:taylorMstar}, we obtain
\begin{equation}
e^{C_n(h)}\Ps\bigl(H_n\ge A_n+x\bigr)
=
\frac{e^{-\tilde\lambda x}+o(1)}{\tilde\lambda},
\qquad x\in B,
\label{eq:tailas}
\end{equation}
uniformly for \(x\) in compact sets. Therefore, for every compact interval
\([a,b]\),
\begin{equation}
\label{eq:kernel_lim}
K_n(h,[a,b])
\xrightarrow{\Ph\text{-a.s.}}
\tilde\lambda^{-1}\bigl(e^{-\tilde\lambda a}-e^{-\tilde\lambda b}\bigr)
=
\int_a^b e^{-\tilde\lambda x}\,dx,
\end{equation}
which is precisely the vague convergence of \(K_n(h,\cdot)\) to
\(\mathcal D_{\tilde\lambda}\).

\section*{REM universality of the energy levels}
We will now restrict the previous kernel asymptotic to the choice
\(C_n(h)=n c_\rho\). Taking the centering \(A_n(h)\) defined by \eqref{eq:Anchoice} for this value of \(C_n(h)\), we prove that the point process \(\mathbf H_n\) converges in distribution to a Poisson point process with intensity measure \eqref{eq:exponential}.

The proof is based on computing the limit Laplace functional of \(\mathbf H_n\). 

Fix a bounded, continuous, nonnegative, compactly supported test function \(f\).
Write \(\Omega_n:=\{-1,1\}^n\). Denoting with $\EU$ the expectation on $U$, we get
\begin{equation}
\EU\!\left[e^{-\langle f,\mathbf H_n\rangle}\right]
=
\prod_{\sigma\in\Omega_n}
\Bigl(1+q_{\sigma,n}\bigl(e^{-f(X_n(\sigma))}-1\bigr)\Bigr),
\label{eq:Laplace-prod}
\end{equation}
where
\begin{equation}
\label{eq:X}
q_{\sigma,n}:=\Qn(\sigma),
\quad
X_n(\sigma):=H_n(h,\sigma)-A_n(h).
\end{equation}
Since \(\Ps(\sigma)\le ((1+|m|)/2)^n\), \eqref{eq:Qn} gives
\begin{equation}
q_{\sigma,n}\le e^{-\delta n},
\quad
\delta=(1-\rho)\bigl(\log 2-\log(1+|m|)\bigr)>0.
\label{eq:qsmall}
\end{equation}
Hence, uniformly for \(z\in[-1,0]\) and \(0\le q\le e^{-\delta n}\),
\[
\log(1+qz)=(1+O(e^{-\delta n}))qz.
\]
Applying this with \(z=e^{-f(X_n(\sigma))}-1\) and summing over \(\sigma\), we
find
\[
\begin{aligned}
\log \EU\!\left[e^{-\langle f,\mathbf H_n\rangle}\right]
&=(1+O(e^{-\delta n}))\sum_{\sigma\in\Omega_n} q_{\sigma,n}\bigl(e^{-f(X_n(\sigma))}-1\bigr)\nonumber\\
&=(1+o(1))\int_{\R}(e^{-f(x)}-1)K_n(h,dx),
\label{eq:Laplace-expansion}
\end{aligned}
\]
Using now the vague convergence of \(K_n(h,\cdot)\), we obtain
\[
\int_{\R}(e^{-f(x)}-1)K_n(h,dx)
\to
\int_{\R}(e^{-f(x)}-1)e^{-\tilde\lambda x}\,dx
\quad \Ph\text{-a.s.}
\]
Since the conditional Laplace functional is bounded by \(1\), dominated
convergence yields
\[
\E\!\left[e^{-\langle f,\mathbf H_n\rangle}\right]
\to
\exp\!\left\{\int_{\R}(e^{-f(x)}-1)e^{-\tilde\lambda x}\,dx\right\}.
\]
This is exactly the Laplace functional of the Poisson point process with
intensity \(\mathcal D_{\tilde\lambda}\), and therefore
\[
\mathbf H_n\xrightarrow{d} \PPP(\mathcal D_{\tilde\lambda}).
\]
\section*{A 1RSB structure for the Gibbs weights}
We now apply the convergence of \(\mathbf H_n\) to determine the limiting joint
distribution of the ranked Gibbs weights in the low-temperature regime
\(\beta>\tilde\lambda\). The proof follows the same strategy as in
\cite[pp.~13--19]{Talcavi} and \cite[Section~4]{ArguinKistler}; see also
\cite[Corollary~{\normalfont Convergence to Poisson--Dirichlet}]{CF26}.

We first truncate the exponential weights, so that the collection of truncated Gibbs weights is a bounded continuous functional of the point
process, and then show that the error due to the truncation is asymptotically negligible.

Let $X_n$ be the recentered Hamiltonian defined in \eqref{eq:X}. We have
\[
G_n(\sigma)=
\frac{e^{\beta X_n(\sigma)}}
{\sum_{\tau\in\mathcal G_n(U)} e^{\beta X_n(\tau)}},
\qquad \sigma\in\mathcal G_n(U).
\]
Let \(w^{(n)}=(w^{(n)}_\alpha)_{\alpha\ge1}\) be the nonincreasing rearrangement
of these Gibbs weights, with $w^{(n)}_\alpha=0$ for $\alpha>|\mathcal G_n(U)|$.

For \(M>0\), define the truncated Gibbs weights
\[
G_n^{(M)}(\sigma)=
\frac{e^{\beta X_n(\sigma)}\1{|X_n(\sigma)|\le M}}
{\sum_{\tau\in\mathcal G_n(U)} e^{\beta X_n(\tau)}\1{|X_n(\tau)|\le M}},
\]
whenever the denominator is nonzero, and \(G_n^{(M)}(\sigma)=0\) otherwise. Define
\[
w^{(n,M)}=(w^{(n,M)}_\alpha)_{\alpha\ge1}
\]
analogously to $w^{(n)}$, from the truncated Gibbs weights.

For fixed \(M\), the map
\(\mathbf H_n\mapsto w^{(n,M)}\) is bounded and continuous at point
measures with no atom at \(\pm M\). Since the limiting Poisson point process
has no atom at \(\pm M\) almost surely, then the continuous mapping theorem gives
\begin{equation}
\label{eq:Gconv}
w^{(n,M)}
\xrightarrow{d}
w^{(M)}:=
\left(
\frac{e^{\beta \eta_\alpha}\1{|\eta_\alpha|\le M}}
{\sum_j e^{\beta \eta_j}\1{|\eta_j|\le M}}
\right)^{\downarrow}_{\alpha\ge1},
\end{equation}
where \(\sum_{\alpha\ge1}\delta_{\eta_\alpha}\) is a Poisson point process with
intensity \(\mathcal D_{\tilde\lambda}\), with the convention that \(w^{(M)}=0\) if the denominator is zero. The symbol $^{\downarrow}$ denotes the rearrangement in decreasing order. 

It remains to remove the cutoff. Set
\[
\begin{aligned}
&S_n:=\sum_{\sigma\in\mathcal G_n(U)}e^{\beta X_n(\sigma)},\\
&S_n(M):=\sum_{\sigma\in\mathcal G_n(U)} e^{\beta X_n(\sigma)}
\times\1{|X_n(\sigma)|\le M},\\
&T_n^-(M):=\sum_{\sigma\in\mathcal G_n(U)} e^{\beta X_n(\sigma)}\times\1{X_n(\sigma)<-M},\\
&T_n^+(M):=\sum_{\sigma\in\mathcal G_n(U)} e^{\beta X_n(\sigma)}\times\1{X_n(\sigma)>M}.
\end{aligned}
\]
For \(\beta>\tilde\lambda\), the left tail is integrable. The sharp
large-deviation estimates used in \cite{CF26} to prove the kernel asymptotics
imply, \(\Ph\)-almost surely,
\[
\lim_{M\to\infty}\limsup_{n\to\infty}
\int_{-\infty}^{-M} e^{\beta x}K_n(h,dx)=0.
\]
Equivalently,
\[
\lim_{M\to\infty}\limsup_{n\to\infty}\EU T_n^-(M)=0.
\]
The right tail is controlled through the number of points. Indeed, by
\eqref{eq:tailas},
\[
\begin{aligned}
\PU\bigl(T_n^+(M)>0\bigr)
&\le \EU\,\mathbf H_n((M,\infty))=K_n(h,(M,\infty))\\
&\xrightarrow{\Ph\text{-a.s.}}
\int_M^\infty e^{-\tilde\lambda x}\,dx
=\frac{e^{-\tilde\lambda M}}{\tilde\lambda}
\xrightarrow[M\to\infty]{}0.
\end{aligned}
\]
Since
\[
S_n=S_n(M)+T_n^-(M)+T_n^+(M),
\]
we conclude that
\[
\lim_{M\to\infty}\limsup_{n\to\infty}
\PU\bigl(|S_n-S_n(M)|>\varepsilon\bigr)=0,
\qquad \varepsilon>0.
\]
Since $S_n>0$, the above inequality implies that, for every \(\delta>0\) one can choose \(M\) and then \(\eta>0\) such that
\[
\liminf_{n\to\infty}\PU\bigl(S_n(M)>\eta\bigr)>1-\delta.
\]
On the event \(\{S_n(M)>\eta\}\), an immediate computation yields
\[
\|w^{(n)}-w^{(n,M)}\|_1
\le 2\,\frac{S_n-S_n(M)}{S_n}.
\]
Therefore
\[
\lim_{M\to\infty}\limsup_{n\to\infty}
\PU\bigl(\|w^{(n)}-w^{(n,M)}\|_1>\varepsilon\bigr)=0,
\qquad \varepsilon>0,
\]
Since also \(w^{(M)}\to w\) almost surely in \(\ell^1\) when \(M\to\infty\), the convergence \eqref{eq:Gconv} extends to $w^{(n)}$:
\[
w^{(n)}
\xrightarrow{d}
w=
\left(
\frac{e^{\beta\eta_\alpha}}{\sum_j e^{\beta\eta_j}}
\right)^{\downarrow}_{\alpha\ge1}.
\]
By \cite[Prop.~10]{PitmanYor} \(w\) has Poisson--Dirichlet law
\(\mathrm{PD}(\tilde\lambda/\beta,0)\). 

The Poisson--Dirichlet limit also yields an explicit asymptotic formula for a
class of Gibbs averages. Specifically, whenever
\((f(\sigma))_{\sigma\in\{-1,1\}^n}\) is a family of positive bounded i.i.d.\
random variables, independent of \(h\) and \(U\), then for
\(\beta>\tilde\lambda\),
\[
\E\!\left[\log\!\left(\sum_{\sigma\in\mathcal G_n(U)}G_n(\sigma)\,f(\sigma)\right)\right]
\to
\frac{\beta}{\tilde\lambda}\log \E\!\left[f(\sigma)^{\tilde\lambda/\beta}\right].
\]
This is the standard \(1\)RSB mechanism of Parisi theory
\cite[Theorem~5.2]{As2}.  Interestingly, the same Gibbs average can also be read as a
single-query cross-attention map with query \(Q=\beta h\), keys
\(K=(\sigma-m\mathbf 1)_{\sigma\in\mathcal G_n(U)}\), and values
\(V=(f(\sigma))_{\sigma\in\mathcal G_n(U)}\) \cite{Aallneed}.

\section*{Free energy}
In this section we will compute the limiting free energy. The resulting formula shows
a phase transition at the critical inverse temperature \(\beta=\tilde\lambda\).
By the previous section, this is precisely the threshold at which the Gibbs
weights enter the Poisson--Dirichlet, or \(1\)RSB, regime.

For \(E\in\R\), let
\[
\mathcal N_n(E):=\mathbf H_n([nE,\infty)).
\]
A standard slicing argument then gives
\begin{equation}
\label{eq:slicing}
\begin{aligned}
F(\beta)&=\lim_{n\to\infty}\left\{\beta\frac{A_n(h)}{n}+\frac1n\log \sum_{\sigma\in\mathcal G_n(U)}
 e^{\beta(H_n(h,\sigma)-A_n(h))}\right\}\\
&\to \beta\tilde a+\sup_{E\in \R}\left\{\beta E+\lim_{n\to \infty}\frac{1}{n}\log \mathcal N_n(E)\right\}.
\end{aligned}
\end{equation}
Conditionally on $h$, $\mathcal N_n(E)$ is a sum of independent Bernoulli variables. Put
\[
\overline{\mathcal N}_n(E):=\EU[\mathcal N_n(E)]=e^{nc_\rho}\Ps\bigl(H_n(h,\sigma)\ge A_n(h)+nE\bigr).
\]
If
\[
p_{\sigma,E}:=q_{\sigma,n}\1{H_n(h,\sigma)\ge A_n(h)+nE},
\]
then, conditionally on $h$,
\[
\mathcal N_n(E)=\sum_{\sigma\in\Omega_n}B_{\sigma,E},
\quad
B_{\sigma,E}\sim\mathrm{Bernoulli}(p_{\sigma,E}),
\]
with the variables $B_{\sigma,E}$ independent. Hence
\begin{equation}
\label{eq:var}
\begin{aligned}
\var_U\bigl(\mathcal N_n(E)\bigr)
&=\sum_{\sigma\in\Omega_n}p_{\sigma,E}(1-p_{\sigma,E})\\
&\le \sum_{\sigma\in\Omega_n}p_{\sigma,E}
=\overline{\mathcal N}_n(E).
\end{aligned}
\end{equation}
The quenched large-deviation principle for \(H_n/n\), together with the limit \eqref{eq:limits}, yields
\begin{equation}
\label{eq:almost_S}
\frac1n\log \overline{\mathcal N}_n(E)\to \sup_{x>\tilde a+E}(c_\rho-G^*(x))=:\mathcal{S}(E),
\,\, \Ph\text{-a.s.}
\end{equation}
where we applied the convention $G^*(x)=\infty$ if $x$ is outside the natural domain of $G^*$. Since $G^*$ attains its minimum at $0$ and it is strictly increasing on the positive branch $(0,\varsigma)$,
\begin{equation}
\label{eq:limits_N}
\mathcal{S}(E)=
\begin{cases}
c_\rho,& \textup{if }E<-\tilde a\\
c_\rho-G^*(\tilde a+E),& \textup{if }E\in [-\tilde a,\varsigma-\tilde a]\\
-\infty,& \textup{if }E>\varsigma-\tilde a
\end{cases}
\end{equation}
If $E>0$, then $G^*(\tilde a+E)>G^*(\tilde a)=c_\rho$, so $\overline{\mathcal N}_n(E)$ decays exponentially to $0$. Hence, Markov's inequality and the Borel--Cantelli lemma give, for every fixed $E>0$,
\[
\mathcal N_n(E)\to 0,\quad a.s.
\]

If $E<0$, then $\frac1n\log \overline{\mathcal N}_n(E)$ converges to a positive quantity. Moreover, by \eqref{eq:var} and Chebyshev's inequality, for every $\varepsilon>0$,
\[
\begin{aligned}
&\PU\Bigl(\bigl|\mathcal N_n(E)-\overline{\mathcal N}_n(E)\bigr|
>\varepsilon \overline{\mathcal N}_n(E)\Bigr)\\
&\qquad\le \frac{1}{\varepsilon^2\overline{\mathcal N}_n(E)}
\sim\frac{e^{-nS(E)+o(n)}}{\varepsilon^2}.
\end{aligned}
\]
By eq. \eqref{eq:limits_N} for $E<0$ the right-hand side is exponentially small in $n$ and summable. Therefore, by the Borel--Cantelli lemma and Chebyshev's inequality,
$\mathcal N_n(E)=\overline{\mathcal N}_n(E)(1+o(1))$ almost surely for every fixed $E<0$. Therefore the entropy profile is
\[
\frac1n\log \mathcal N_n(E)\to
\begin{cases}
\mathcal{S}(E), &  E<0,\\[1mm]
-\infty, & E>0,
\end{cases}
\quad a.s.,
\]
where the value at $E=0$ is understood by taking the limit $E\uparrow0$ in the slicing variational formula.
The slicing argument in \eqref{eq:slicing}, formula \eqref{eq:limits_N}, and the above formula give
\[
F(\beta)=
\begin{cases}
c_\rho+G(\beta), & \beta\le \tilde\lambda,\\[1mm]
\beta\tilde a, & \beta>\tilde\lambda,
\end{cases}
\qquad \Ph\text{-a.s. and }\PU\text{-a.s.}
\]
Notice that the second derivative is discontinuous at \(\beta=\tilde\lambda\). Thus the model
undergoes the phase transition at this critical inverse temperature. By the
previous section, this is exactly the onset of the Poisson--Dirichlet distribution, or
\(1\)RSB, behavior of the Gibbs weights. 

\section*{Conclusion}

We have shown that a purely linear random Hamiltonian exhibits REM universality after an independent random selection that retains an exponential number of configurations. After an environment-dependent centering, the selected energy levels converge to a Poisson point process with exponential intensity. As a consequence, in the low-temperature regime the ranked Gibbs weights converge to the Poisson--Dirichlet law \(\mathrm{PD}(\tilde\lambda/\beta,0)\), equivalently to a one-level \(1\)RSB Ruelle probability cascade, and the free energy is explicit, with a freezing transition at \(\beta=\tilde\lambda\). These results extend previous forms of REM universality beyond shrinking energy windows, sparse random selection, and models with an added REM term. They show that REM statistics, Poisson--Dirichlet Gibbs structure, and the associated thermodynamic transition already emerge in a genuinely correlated linear model. Full proofs and finite-\(n\) estimates are given in \cite{CF26}.

\end{document}